\documentclass[12pt]{article}
\usepackage{amsmath}
\usepackage{graphicx}
\usepackage{comment} %comment blocks
\usepackage{xcolor}
\usepackage{amssymb}
\usepackage{cite}
\usepackage{amsthm}
\newcommand{\blind}{0}

\begin{document}

\def\spacingset#1{\renewcommand{\baselinestretch}%
{#1}\small\normalsize} \spacingset{1}

%%%%%%%%%%%%%%%%%%%%%%%%%%%%%%%%%%%%%%%%%%%%%%%%%%%%%%%%%%%%%%%%%%%%%%%%%%%%%%

\if0\blind
{
  \title{\bf The inter-play between asymmetric noise and coupling in a super-critical Hopf bifurcation}
  \author{Gurpreet Jagdev, Na Yu\\
    {\small Department of Mathematics, Toronto Metropolitan University, Toronto, Canada}\\}
\date{}    
  \maketitle
} \fi

\if1\blind
{
  \bigskip
  \bigskip
  \bigskip
  \begin{center}
    {\LARGE\bf Title}
\end{center}
  \medskip
} \fi

%\bigskip
%\begin{abstract}
%[ABSTRACT]
%\end{abstract}

%\noindent%
%{\it Keywords: [ . . . ]}  %4 to 6 keywords, that do not appear in the title , 

%\newpage
\spacingset{1.8} % DON'T change the spacing!

%%%%%%%%%%%%%%%%%%%%%%%%%%%%%%%%%%%%%%%%%%%%%%%%%%%%%%%%%%%%%%%%%%%%%%%%%%%%%%%%%%%%%%%%%%%%%%%%%%%%%%%%%%%%%%%%%%%%%%%%%%%%%%%%%%%%%%%%%%%%%%%%%%%%%%%%%%%%%%%%%%%%%%%%%%%%%%%%%%%%%%%%%%%%%%%%%%%%%%%%%%%%%%%%%%%%%%
\section{Introduction}
\label{sec:intro}
%----------------------
Coupled oscillators are ubiquitous in neuroscience and can be used to model a wide variety of oscillatory phenomena (e.g. modeling a neural pacemaker \cite{schafer1998heartbeat}). Although the dynamics of practical systems rarely, if ever, exhibit perfect periodicity, the dynamics of such systems can nonetheless be understood as the manifestation of a stochastic limit cycle \cite{freund2003frequency,anishchenko2002chaotic}. \par

When multiple oscillators interact, synchronization is often an inescapable phenomenon. In the most general sense, synchronization can be described as the mutual adjustment of the oscillatory rhythms (i.e. phase) of oscillators due to relatively weak interactions. This phenomenon was first discovered by Huygens in the seventeenth-century when studying the simple motion of two pendulum clocks \cite{oliveira2015huygens}. Huygens concluded that the apparent synchronization was likely caused by the weak interactions between the two pendulum clocks \cite{oliveira2015huygens}. However, when studying the motion of stochastic oscillators, synchronization can also arise as a consequence of the influence of common noise on the system of oscillators  by means of a resonant type mechanism \cite{garcia2009competition}.\par

In this chapter we study the interplay between noise and coupling in a super-critical HB. We consider a pair of diffusively coupled $\lambda-\omega$ oscillators with parameters chosen such that the model is in the vicinity of a super-critical HB, quiescent in the absence in noise, and excitable with the addition of an intrinsic noise stimulus. Our results agree with previous studies (e.g. \cite{yu2006stochastic,thompson2012stochastic,yu2021noise,yu2008stochastic}) which show that noise can play a constructive role in the PS of coupled oscillators. Many such studies tend to assume symmetrical interactions between pairs of oscillators (i.e. symmetrical coupling and/or stochastic stimulus). However, in biological systems the interactions between such oscillators are often asymmetric and the assumption of symmetric interactions is quite restrictive (e.g. see \cite{sheeba2009asymmetry,cimponeriu2003inferring}). To address this issue we allow both coupling and noise to be asymmetric. We find that the asymmetries between the couplings and noise have a robust effect on PS and remarkably, that symmetrical coupling and noise lead to relatively low levels of PS---all else equal. 
\par 
The remainder of this chapter is organized as follows: Sec. \ref{sec:model} describes the model; Sec. \ref{sec:methods} describes analytic and numerical methods. Sections \ref{noise-induced oscillations} and \ref{sec: CR} study the effects of additive noise on the dynamics and PS of our system. Sections \ref{sec:lambda} and \ref{The effects of coupling on synchronization} discuss the effects of the bifurcation parameter, $\lambda_0$, and coupling strengths, $d_1$ and $d_2$ on PS. A discussion is given in Sec. \ref{sec:discussion}.

%%%%%%%%%%%%%%%%%%%%%%%%%%%%%%%%%%%%%%%%%%%%%%%%%%%%%%%%%%%%%%%%%%%%%%%%%%%%%%%%%%%%%%%%%%%%%%%%%%%%%%%%%%%%%%%%%%%%%%%%%%%%%%%%%%%%%%%%%%%%%%%%%%%%%%%%%%%%%%%%%%%%%%%%%%%%%%%%%%%%%%%%%%%%%%%%%%%%%%%%%%%%%%%%%%%%%%
\section{Model and methods}
\label{sec:model and methods}
\subsection{Model}
\label{sec:model}
We adapt the canonical model for a HB, $\lambda-\omega$ system, to develop a pair of coupled oscillators with additive noise which are modelled by the set of stochastic differential equations (SDEs)
\begin{align}
    dx_i &= [\lambda(r_i)x_i - \omega(r_i)y_i + d_i(x_j-x_i)]dt+\delta_i d\eta_i(t), \label{x_i}\\
    dy_i & = [\omega(r_i)x_i + \lambda(r_i)y_i + d_i(y_j-y_i)]dt, \label{y_i}\\
    r_i^2 &= x_i^2 + y_i^2, \label{r_i}
\end{align}
where $i=1,2$ and $j=2,1$. $r_i = \sqrt{x_i^2 + y_i^2}$ represents the amplitude of the $i$th oscillator. $\lambda(r_i) = \lambda_0 + \alpha r_i^2 + \gamma r_i^4$, controls the increment and decrement of the amplitude of the $i$th oscillator. In particular, $\lambda_0$ is the control parameter and a HB occurs at $\lambda_0=0$; $\alpha$ and $\gamma$ influence the system away from the bifurcation point. $\omega(r_i) = \omega_0 + \omega_1 r_i^2$  determines the increment and decrement of the frequency of the $i$th oscillator, where $\omega_1$ governs the evolution of the frequency with respect to the amplitude, $r_i$. Particularly, the amplitude does not directly effect phase when $\omega_1=0$. We consider a super-critical HB, and accordingly, we choose parameters values: $\alpha=-0.2$; $\gamma=-0.2$; $\omega_0=2$; and $\omega_1=0$ \cite{yu2006stochastic}. The terms $d_i(x_j-x_i)$ and $d_i(y_j-y_i)$ are diffusive coupling terms with coupling strength $d_i$, where, $i=1,2$ and $j=2,1$. The term $\delta_i d\eta_i(t)$ represents an intrinsic noise applied to $x_i$ (i.e. the noise is unique to each oscillator) where the function $\eta_i(t)$ is a Wiener process with zero mean and unity variance (i.e. a standard Brownian motion) and $\delta_i$ is a scaling parameter, also known as noise intensity. To study the dynamics of coupled oscillations we restrict our attention to excitatory coupling, $0.01 \leq d_1,d_2 \leq 0.3$ \cite{yu2008stochastic}.  
%with $r_i=0$ as a fixed point of the deterministic system (i.e. $\delta_1=\delta_2=0$). 
%We study the model presented in Eqns. \ref{x_i} - \ref{r_i} as it is the canonical model for a normal for near a Hopf bifurcation (HB). To ensure that our model undergoes a super-critical HB as $\lambda_0$ traverses the bifurcation point in the deterministic regime ($\delta_1=\delta_2=0$) 

%--------------------------------------------------------------------------------------------------
\subsection{Methods}
\label{sec:methods}
 To study the interplay of noise and coupling in the synchronization of our model we analyze the PS of both oscillators when subject to the additive noise, $\delta_i d\eta_i$, $i=1,2$. Since the oscillators rotate about the fixed point $(x_i,y_i)=(0,0)$, $i=1,2$, when driven by noise, the phase of $x_i$ is taken to be the natural phase \cite{rosenblum2001phase},
 \begin{equation}
 \label{naturalphase}
      \phi_i = \arctan{(y_i/x_i)},
  \end{equation}
$i=1,2$. In the classical treatment of phase analysis, PS measures are often based on the distribution of the phase difference, $\Delta \phi_{n,m} = n\phi_2-m\phi_1$, where $n,m \in \mathbb{N}$ characterize the order of locking \cite{rosenblum2001phase, rosenblum2022synchronization}. However, in the presence of noise the phase of the oscillators can exhibit random jumps of $\pm 2\pi$, called phase slips, which can cause the phase difference, $\Delta \phi_{n,m}$, to become unbounded and moreover, lead to inaccurate numerical results. Therefore, instead of considering the natural phase in Eqn. \ref{naturalphase}, we consider the cyclic relative phase (CRP) \cite{rosenblum2001phase,pikovsky2002synchronization, zaks1999alternating,mormann2000mean}, 
 \begin{equation}
 \label{phase}
     \varphi_i = \phi_i\text{ mod } 2\pi,
 \end{equation}
 which is the natural phase wrapped $2\pi$. This ensures that the  phase difference, $\Delta \varphi_{n,m}$, will be bounded. For simplicity, we consider only $1-1$ synchronization, that is, $\Delta \varphi=\Delta \varphi_{1,1}$. \par
 %%%%%%%%%%%%%%%%%%%%%%%%%%%%%%%%%%%%%%%%%%
 
The bifurcation diagrams (i.e. Fig. 3.1 and 3.2) are generated using XPPAUT software \cite{ermentrout2012xppaut}. All further analysis (i.e. Fig. 3.3 - 3.9) is conducted using MATLAB. To simulate the SDEs in Eqns. \ref{x_i} - \ref{r_i} we use the Euler-Maruyama method over the time range $[15, 100]$ with time-step $dt=0.01$ and arbitrary random initial conditions $x_i(0),y_i(0) \sim N(0,0.008^2)$, $i=1,2$. Finally, high-frequency fluctuations are removed from the time series of $x_1$ and $x_2$ by applying a low-pass filter. The signal-to-noise ratio, $\beta$, and synchronization measures $|\Delta \varphi|$, $R$, and $\rho$ in Eqns. \ref{beta} - \ref{rho}, respectively, are averaged over $N=200$ simulations and time domain $t \in [15,100]$ (in arbitrary units).

%In order to remove high-frequency fluctuations from the time series of $x_1$ and $x_2$ (which results in more consistent numerical results) we smooth $x_1$ and $x_2$ using local regression with weighted linear least squares and a $2$nd degree polynomial model using the built-in MATLAB function \textit{smooth($\sim$,0.003,'loess')}. All, further analysis--- such as the computation of phase, synchronization/coherence measures, etc.--- is done using MATLAB. The signal-to-noise ratio $\beta$ in Eqn. \ref{beta}, and synchronization measures $|\Delta \varphi|$, $R$, and $\rho$ in Eqns. \ref{absolute phase difference}, \ref{R}, and \ref{rho}, respectively, are averaged over $N=200$ simulations and time domain $15 \leq t \leq 100$. Note, the component $t \leq 15$ is removed to ensure the system has settled to an equilibrium in our analysis. Lastly, to compute the synchronization measure $\rho$, the MATLAB function \textit{wentropy($\sim$,'shannon')} is used to compute $S$ and $S_{max}$ in Eqn. \ref{rho}. \par
%----------------
%%%%%%%%%%%%%%%%%%%%%%%%%%%%%%%%%%%%%%%%%%%%%%%%%%%%%%%%%%%%%%%%%%%%%%%%%%%%%%%%%%%%%%%%%%%%%%%%%%%%%%%%%%%%%%%%%%%%%%%%%%%%%%%%%%%%%%%%%%%%%%%%%%%%%%%%%%%%%%%%%%%%%%%%%%%%%%%%%%%%%%%%%%%%%%%%%%%%%%%%%%%%%%%%%%%%%%
\section{Results}

\subsection{Bifurcation analysis}
\label{bifurcation structure}

We consider three cases of the deterministic system ($\delta_1=\delta_2=0$): single oscillators (i.e. uncoupled oscillators) with $d_1=d_2=0$ (Fig. 3.1a); two symmetrically coupled oscillators with $d_1=d_2\neq 0$ (Fig. 3.1b); and two asymmetrically coupled oscillators with $d_1 \neq d_2$ (Fig. 3.1c and 3.1d). When $\lambda_0<0$, both oscillators in all cases are quiescent (i.e. stable fixed points) in the absence of noise (black solid line in Fig. 3.1). Conversely, when $\lambda_0>0$, stable periodic orbits emerge (blue solid line), and $x_1=x_2=0$ are unstable fixed points (black dashed line). Hence, both the single (Fig. 3.1a) and coupled (Fig. 3.1b - 3.1d) oscillators undergo a supercritical HB (denoted as $\text{HB}_1$) at $\lambda_0=0$. However, when the oscillators are coupled the system exhibits a second HB (denoted as $\text{HB}_2$) which leads to unstable periodic orbits (blue dashed line in Fig. 3.1b - 3.1d). 
%: the first occurs at $\lambda_0=0$ which results in a stable periodic orbit; and the second occurs at $\lambda_0 \approx 0.1$ which results in an unstable periodic orbit (see the dashed blue branches in Fig. 1b and 1c). 
When the oscillators are symmetrically coupled (e.g. Fig. 3.1b) the amplitude of the periodic orbits generated by both oscillators are identical. The unstable periodic orbit has an amplitude which is slightly less than the amplitude of the stable periodic obit, and with the increment of $\lambda_0$  both converge. When the oscillators are coupled asymmetrically (e.g. $d_1=0.1$ and $d_2=0.01$ in Fig. 3.1c), the bifurcation diagram for oscillator $x_1$ behaves the same as Fig. 3.1b. On the other hand, the amplitude of the unstable periodic orbit for the oscillator $x_2$ (Fig. 3.1d) is near zero.\par

\begin{figure}
\begin{center}
\includegraphics[width=6in]{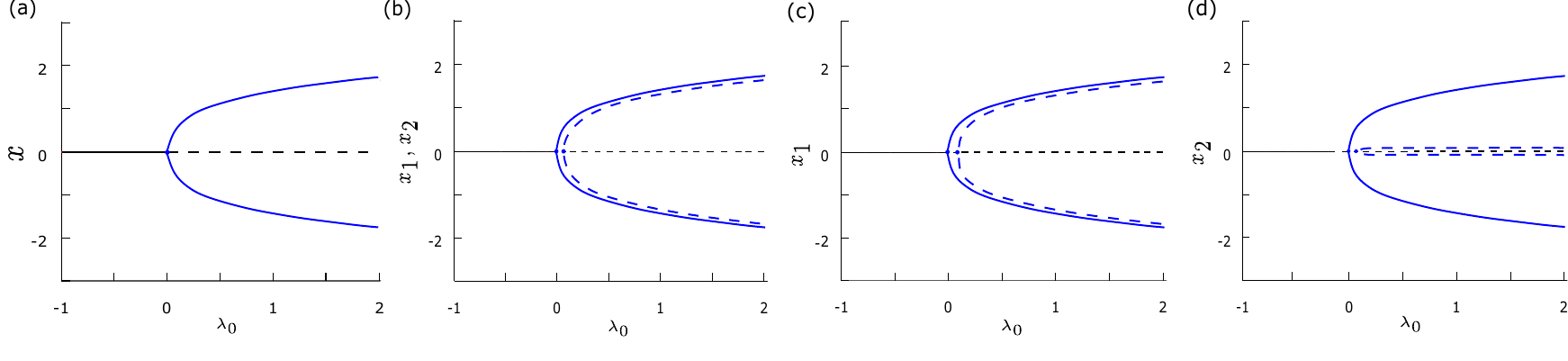}
\end{center}
\caption{Bifurcation diagrams vs. $\lambda_0$ for: (a) single/uncoupled oscillator ($d_1=d_2=0$) in the deterministic regime ($\delta_1=\delta_2=0$); (b) symmetrically coupled oscillators ($d_1=d_2=0.05$) in the deterministic regime; (c) asymmetrically coupled oscillator, $x_1$ vs. $\lambda_0$ ($d_1=0.1$ and $d_2=0.01$), in the deterministic regime; and (d) asymmetrically coupled oscillator, $x_2$ vs. $\lambda_0$ ($d_1=0.1$ and $d_2=0.01$), in the deterministic regime. Note that panel c shows that bifurcation diagram for oscillator $x_2$. Stable solutions are marked by solid blue lines and solid black lines, unstable solutions are marked by dashed blue lines and dashed black lines. Other parameters are: $\alpha=-0.2$; $\gamma=-0.2$; $\omega_0=2$; and $\omega_1=0$.}
\end{figure}

To further explore the effect(s) of coupling on our model (in the deterministic regime; $\delta_1=\delta_2=0$) we calculated the two-parameter bifurcation diagrams for $x_1$ taking the coupling strengths $d_1$ and/or $d_2$ and $\lambda_0$ as the control parameters. The diagrams are presented in Fig. 3.2, where each point on a line marks a double, $\{\lambda_0,d_i\}$, or triple, $\{\lambda_0,d_1,d_2\}$, of parameters which correspond to a critical point of our model. Recall, when the oscillators are coupled, the model undergoes two HBs (e.g. Fig. 3.1). Accordingly, in Fig. 3.2, we label the points corresponding to the first HB as $\text{HB}_1$ and points corresponding to the second HB as $\text{HB}_2$. Let us first consider the behaviour of the second critical point which corresponds to the branch $\text{HB}_2$. From the two-parameter bifurcation diagram of the symmetrically coupled oscillators shown in Fig. 3.2a, we see that as the coupling strength is increased, the critical points on $\text{HB}_2$ increase with the common coupling strength $d=d_1=d_2$. Indeed, the slope of $\text{HB}_2$ is approximately $0.5$ which indicates that the second Hopf point occurs at twice the coupling strength, $\lambda_0 = 2d$. The bifurcation diagrams of the asymmetrically coupled oscillators are shown in Fig. 3.2b and 3.2c. In Fig. 3.2b, we fix $d_2=0.05$ and vary $d_1$ and in Fig. 3.2c we fix $d_1=0.05$ and vary $d_2$. Considering the branches $\text{HB}_2$ in both panels b and c, we see that both have a slope of one. This indicates that the second Hopf point shifts with the increment of $d_1$ and $d_2$ such that $\lambda_0=d_1$ and $\lambda_0=d_2$, respectively, are critical points. Comparing the two modes of coupling, we see that when the coupling is symmetric and $d=d_1=d_2$ increases, both $d_1$ and $d_2$ increase and thus the critical point is shifted to $\lambda_0=d_1+d_2=2d$. Conversely, when the coupling is asymmetric and either $d_1$ or $d_2$ are increased, either $d_2$ or $d_1$ must remain fixed, and thus the critical point shifts to $\lambda_0=d_1 \text{ or } d_2$.\par

Finally, let us consider the branch $\text{HB}_1$ in Fig. 3.2 (see solid blue branch in Fig. 3.2). We see that for all coupling regimes (Fig. 3.2a - 3.2c), there is always a critical point at $\lambda_0=0$. Furthermore, it follows that whenever $\lambda_0<0$ our system has the single fixed point $x_1,x_2=0$ for any mode of coupling and undergoes a super-critical HB as $\lambda_0$ traverses  the critical point $\lambda_0=0$.\par

\begin{figure}
\begin{center}
\includegraphics[width=6in]{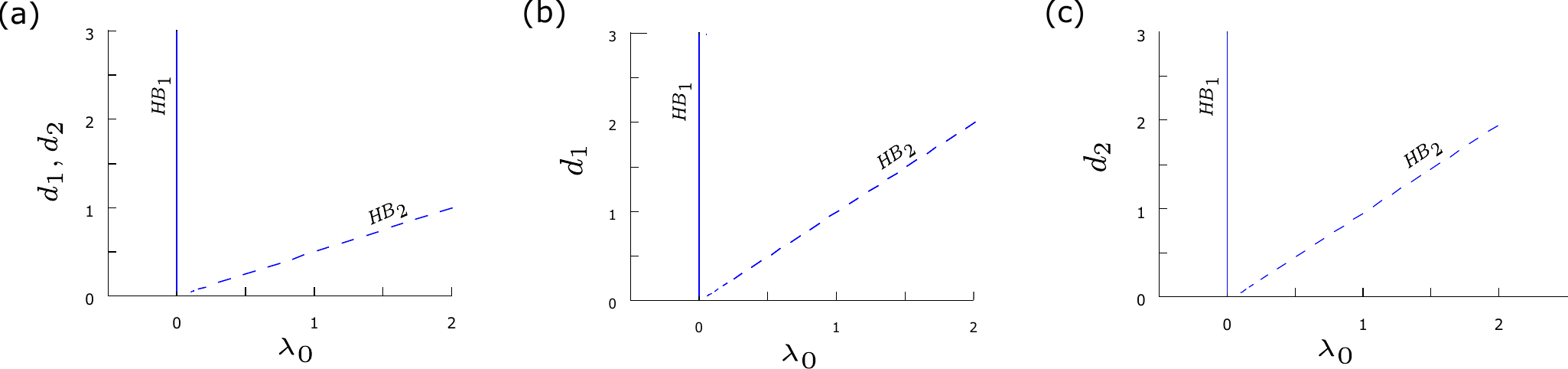}
\end{center}
\caption{Two-parameter bifurcation diagrams for $x_1$ under various coupling regimes: (a) symmetric coupling, $d_1=d_2$ vs. $\lambda_0$; (b) asymmetric coupling, $d_1$ vs. $\lambda_0$ with $d_2=0.05$; and (c) asymmetric coupling, $d_2$ vs. $\lambda_0$ with $d_1=0.05$.  The branches labelled $\text{HB}_1$ (solid blue) and $\text{HB}_2$ (dashed blue) correspond to two distinct HBs for the coupled $\lambda-\omega$ oscillator $x_1$. Other parameters are: $\delta_1=\delta_2=0$; $\alpha=-0.2$; $\gamma=-0.2$; $\omega_0=2$; and $\omega_1=0$.}
\end{figure}
%--------------------------------------------------------------------------------------------------

\subsection{Noise-induced oscillations}
\label{noise-induced oscillations}
%We study the diffusively coupled $\lambda-\omega$ system in Eqns. \ref{x_i} - \ref{r_i} with the parameters listed in Sec. \ref{sec:model} and $\lambda_0=-0.5$. As shown section \ref{sec:model}, when $\lambda_0<0$ our system is quiescent in the absence of noise. 
In this section, we study the noise-induced oscillations of our model. Hence, the control parameter $\lambda_0$ must be in the excitable regime; we let $\lambda_0=-0.5$. The deterministic system (i.e. $\delta_1=\delta_2=0$ as in Fig. 3.3a ) exhibits damped oscillations which converge to the fixed point $x_1=x_2=0$. With the addition of the intrinsic noise, $\delta_i d\eta_i$, $i=1,2$, random perturbations can cause excursions from the stable fixed point $x_1=x_2=0$, which results in oscillatory motion (i.e. noise-induced oscillations) by means of a CR type mechanism \cite{yu2006stochastic,yu2008stochastic,thompson2012stochastic}. Examples can be seen from Fig. 3.3b - 3.3d which show sample times series of $x_1$ and $x_2$ in the presence of noise with different intensities. \par

When the noise intensity is small and symmetric (e.g. $\delta_1=\delta_2=0.01$ in Fig. 3.3b) there are intermittent periods of phase drift and phase locking (both in-phase and anti-phase). When noise intensities are made asymmetric by increasing one of noise intensities, for example, $\delta_1=0.01$ and $\delta_2=0.05$ as in Fig. 3.3c, the time series for $x_1$ and $x_2$ exhibit increased regularity and appear to be better in-phase. For example, there are longer epochs of phase locking. However, if $\delta_2$ is increased further, for example, $\delta_2=3$ as in Fig. 3.3d, the oscillations of $x_2$ become less regular and more chaotic. Moreover, as $\delta_2$ is increased from an optimal level, PS is reduced. This indicates that PS can be optimized by the tuning the noise intensities $\delta_1$ and $\delta_2$. Additionally, one sees that PS is optimized when the noise levels $\delta_1$ and $\delta_2$ are asymmetric (i.e. $\delta_1/\delta_2 \neq 1$). An example of this may be seen from Fig. 3.3, where $x_1$ and $x_2$ appear to be better in-phase when $\delta_1=0.01$ and $\delta_2=0.05$ (Fig. 3.3c) relative to $\delta_1=\delta_2=0.01$ (Fig. 3.3b).
 %--------------------------------------------------------------------------------------------------
\begin{figure}
\begin{center}
\includegraphics[width=5.8in]{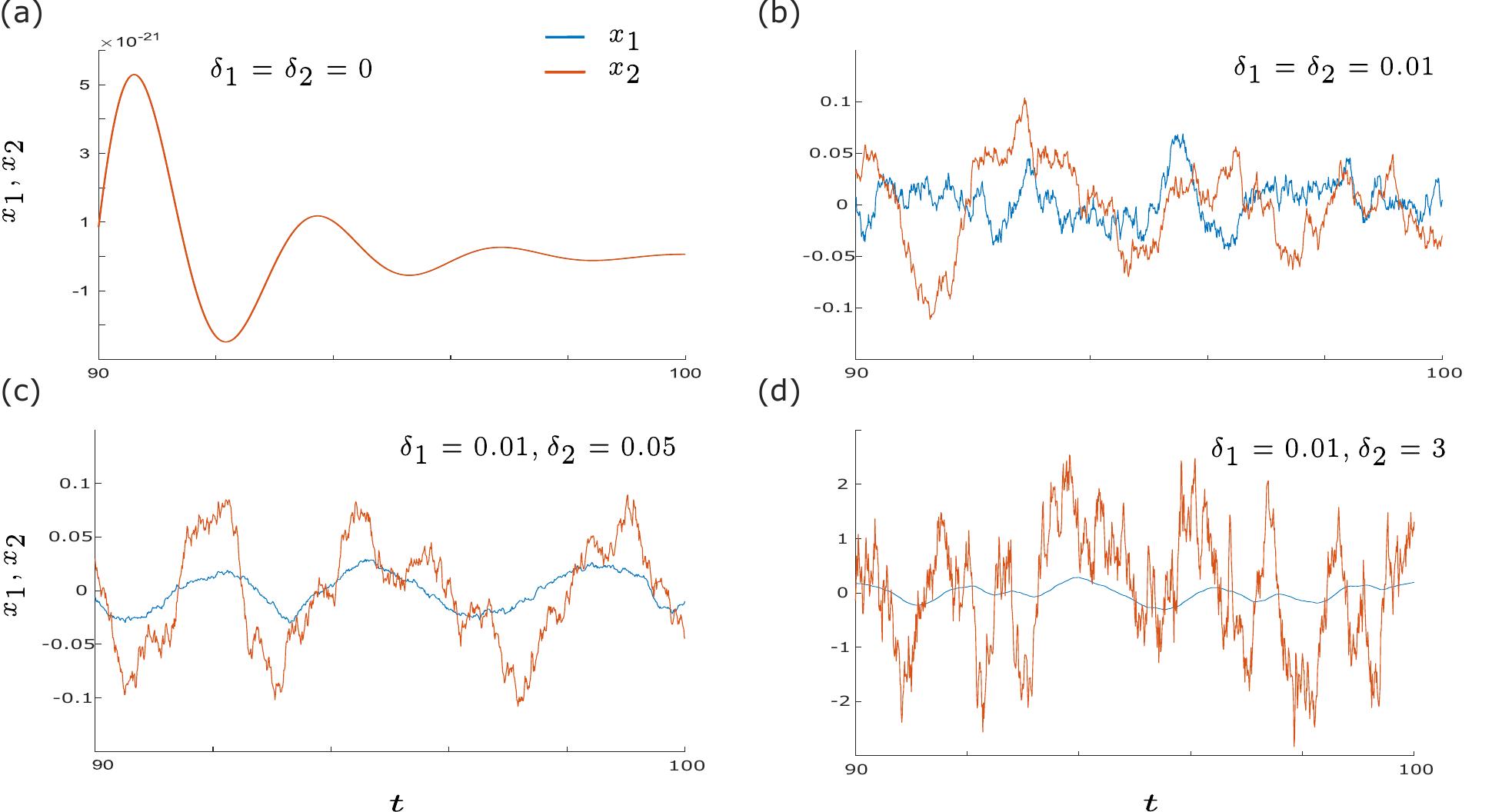}
\end{center}
\caption{Time series of oscillators $x_1 \text{ and } x_2$ under different noise regimes, $x_1,x_2$ vs. $t$. For panel: (a) $\delta_1=\delta_2=0$; (b) $\delta_1=\delta_2=0.01$; (c) $\delta_1=0.01$ and $\delta_2=0.05$; and (d) $\delta_1=0.01$ and $\delta_2=3$. The blue lines represent $x_1$ and the orange lines represent $x_2$. Other parameters are: $\alpha=-0.2$; $\gamma=-0.2$; $\omega_0=2$; $\omega_1=0$; $\lambda_0=-0.5$; $d_1=0.3$; and $d_2=0.01$.}
\end{figure}
%--------------------------------------------------------------------------------------------------
\par
To determine an appropriate range of noise intensities, $\delta_1$ and $\delta_2$, we make use of the signal-to-noise ratio (SNR) measure \cite{gang1993, pikovsky1997},
\begin{equation}
    \beta= h_p(\Delta \omega/\omega_p)^{-1}, \label{beta}
\end{equation} 
where $h_p$ and $\omega_p$ denote the height and central frequency of the power spectrum density (PSD) peak of the osciallator $x$, respectively, and $\Delta\omega$ denotes the width of the PSD peak at half-maximal power, $e^{-\frac{1}{2}}h$. We compute $\beta$ vs. $\delta$ for a single $\lambda-\omega$ oscillator, $x$, and present the results in Fig. 3.4. When $\delta<0.01$, $|\beta| << 1$, which indicates that the level of noise is too weak. As the level of noise is increased the $\beta$ curve exhibits a peak at $\delta \approx 1$ which indicates that this is the optimal noise intensity. Conversely, for the upper range of $\delta$, i.e. $\delta > 1$, $\beta$ is sharply decreasing, which indicates that the noise intensity is overpowering the regularity of the oscillators. Accordingly, we consider $0.01\leq \delta_1,\delta_2\leq 5$.\par
%--------------------------------------------------------------------------------------------------
\begin{figure}
\begin{center}
\includegraphics[width=2.8in]{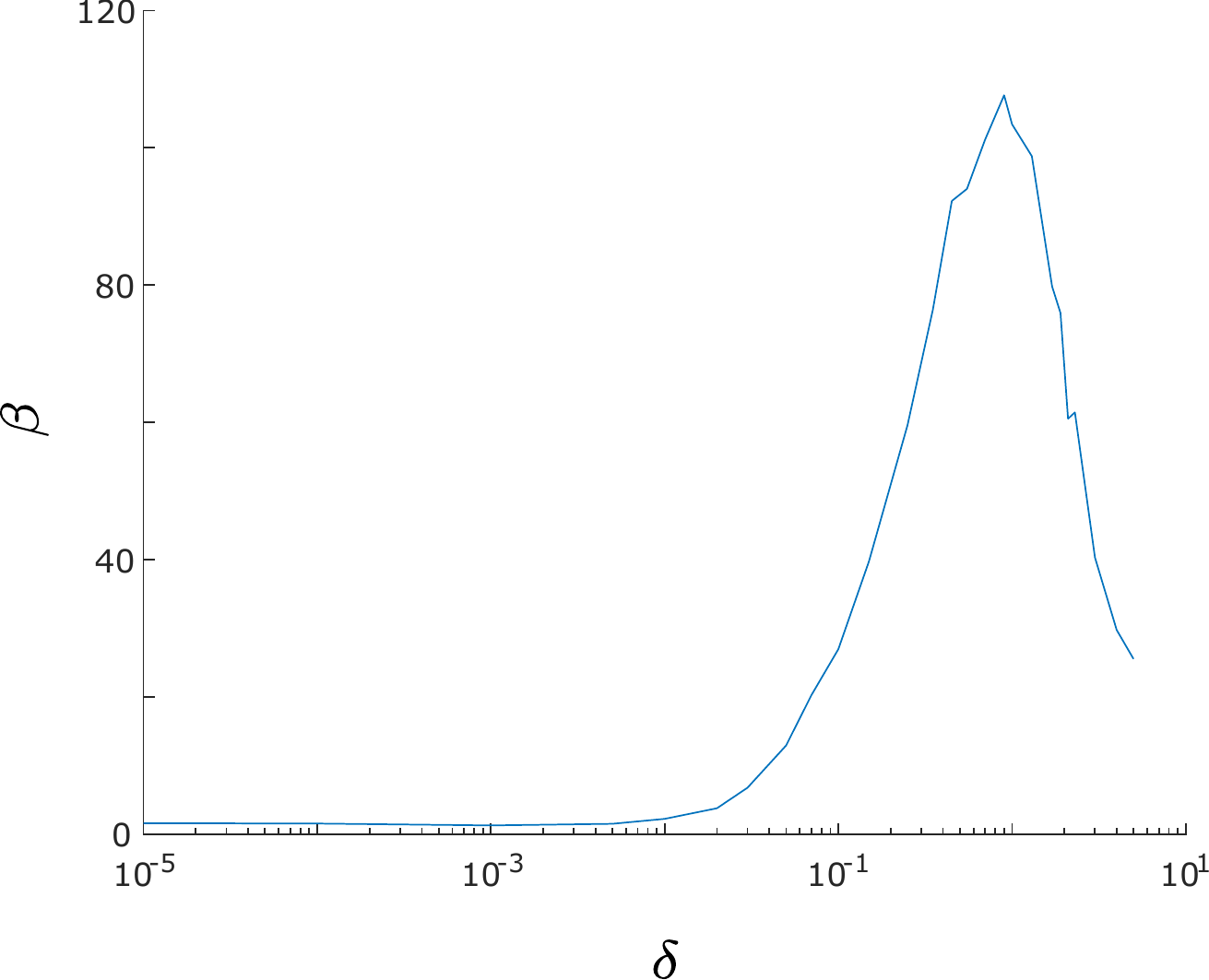}
\end{center}
\caption{SNR measure $\beta$ vs. $\delta$ for a single $\lambda-\omega$ oscillator. Other parameters are: $\alpha=-0.2$; $\gamma=-0.2$; $\omega_0=2$; $\omega_1=0$; and $\lambda_0=-0.5$.}
\end{figure}
%-------------------------------------------------------------------------------------------------------------
%%%%%%%%%%%%%%%%%%%%%%%%%%%%%%%%%%%%%%%%%%%%%%%%%%%%%%%%%%%%%%%%%%%%%%%%%%%%%%%%%%%%%%%%%%%%%%%%%%%%%%%%%%%%%%%%%%%%%%%%%%%%%%%%%%%%%%%%%%%%%%%%%%%%%%%%%%%%%%%%%%%%%%%%%%%%%%%%%%%%%%%%%%%%%%%%%%%%%%%%%%%%%%%%%%%%%%
\subsection{Noise-induced phase synchronization} 
\label{sec: CR}
 The results in Sec. \ref{noise-induced oscillations} indicate that the PS of our model can be optimized by tuning the noise intensities $\delta_1$ and $\delta_2$. To explore the effects of noise on the PS of our model more concretely, we introduce three (time-averaged) measures of PS. First, the  absolute CRP difference, $\Delta \varphi$, defined as
 \begin{equation}
 \label{absolute phase difference}
   \Delta \varphi = \frac{1}{T} \sum^{T}_{t = t_0}|\varphi_1 - \varphi_2|,
 \end{equation}
where $\varphi_i$ is the CRP of the $i$th oscillator. Since we consider excitatory coupling, we study the dynamics of in-phase of oscillators, thus, smaller values of $\Delta \varphi$ relate a greater degree of PS. The second measure we use is the mean phase coherence, $R$, defined as \cite{mormann2000mean, rosenblum2001phase}
\begin{equation}
    \label{R}
    R= \sqrt{\left(\frac{1}{T} \sum^{T}_{t = t_o} \sin{\Delta \varphi}\right)^2 + \left(\frac{1}{T} \sum^{T}_{t = t_o} \cos{\Delta \varphi}\right)^2},
\end{equation}
where $t_0=15$, $T=100$, and $\Delta\varphi=\varphi_1-\varphi_2$. From Eqn. \ref{R}, one sees that $0 \leq R \leq 1$, and greater values of $R$ indicate a greater degree of PS.  The third synchronization measure we use is the normalized synchronization index, $\rho$, defined as \cite{rosenblum2001phase}
\begin{equation}
\label{rho}
    \rho = \frac{S_{max}-S}{S_{max}},
\end{equation}
where $S= -\sum_{k=1}^N p_k \ln{p_k}$ is the Shannon entropy, $S_{max}=\ln{N}$ is the maximum entropy, $N$ is the number of bins, and $p_k$ is the probability of finding $\Delta \varphi$ in the $k$th bin. $\rho$ is normalized such that $0 \leq \rho \leq 1$, and because $S$ is a measure of entropy, it follows that lower values of $\rho$ correspond to a narrower distribution of $\Delta \varphi$ and therefore a greater degree of PS. \par %Although all of the measures listed above measure PS, they are qualitatively different. $|\Delta \varphi|$ in Eqn. \ref{absolute phase difference} measures the peak of the distribution of $\Delta \varphi$, whereas, $R$ and $\rho$ in Eqn. \ref{R} and Eqn. \ref{rho}, respectively, measure the dispersion of $\Delta \varphi$ around its peak \cite{rosenblum2001phase}.  \par

%--------------------------------------------------------------------------------------------------
\begin{figure}
\begin{center}
\includegraphics[width=6in]{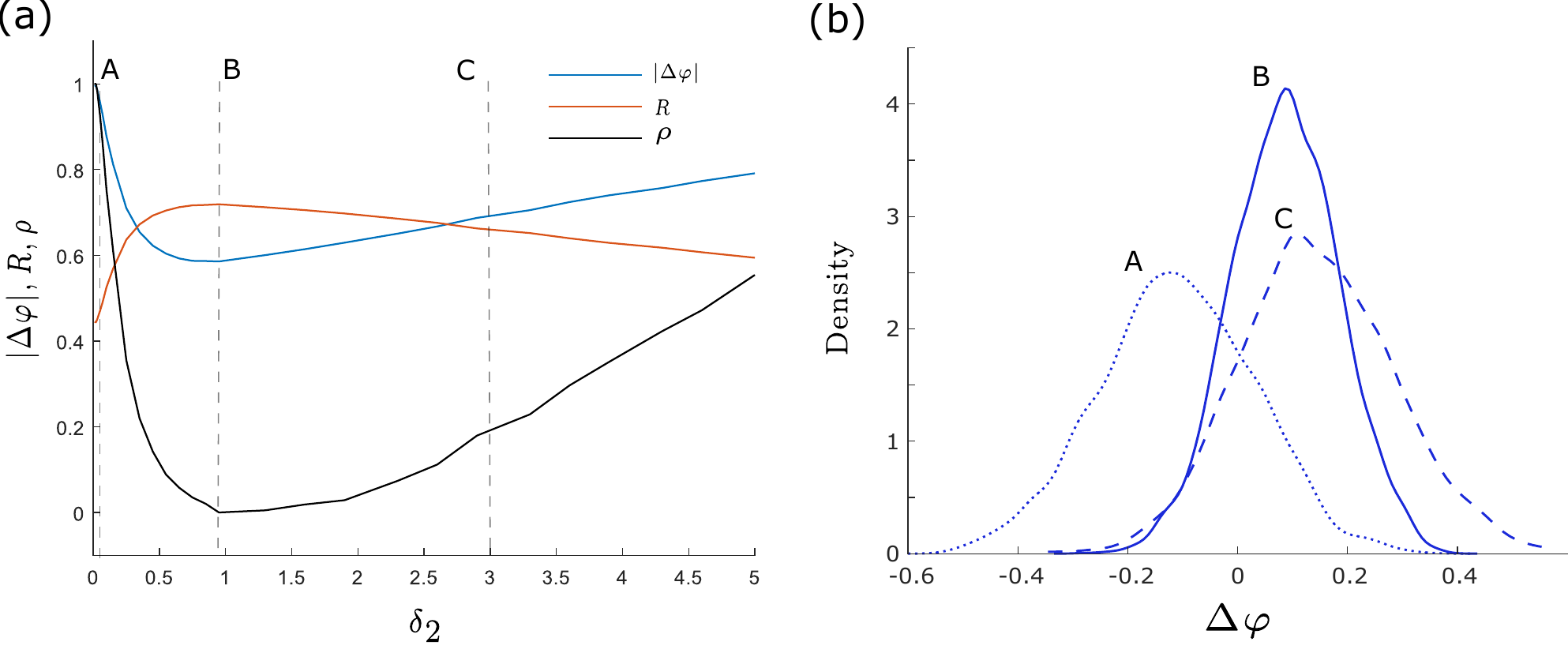}
\caption{(a): absolute CRP difference,  $\Delta \varphi$, mean phase coherence, $R$, and the synchronization index, $\rho$, vs. $\delta_2$. The orange, blue and black curves correspond $R$, $|\Delta \varphi|$, and $\rho$, respectively. The points A, B, and C correspond to $\delta_2 = 0.05, 0.95, \text{ and } 3$, respectively. (b): empirical probability density functions of $\Delta \varphi$. The curves A, B, and C correspond to noise intensities $\delta_2 = 0.05, 0.95, \text{ and } 3$, respectively. $|\Delta \varphi|$, $R$, and $\rho$ are averaged over $200$ trials and empirical probability density functions are averaged of $10$ trials. Other parameters are: $\delta_1=0.05$, $\alpha=-0.2$; $\gamma=-0.2$; $\omega_0=2$; $\omega_1=0$; $\lambda_0=-0.5$; $d_2=0.01$; and $d_1=0.3$.}
\end{center} 
\end{figure}
%--------------------------------------------------------------------------------------------------
To begin our analysis, we fix $\delta_1$ at some (appropriate) arbitrary value (e.g. $\delta_1=0.05$ in Fig. 3.5), simulate our model for $0.01 \leq \delta_2 \leq 5$ with $\lambda_0=-0.5$ and coupling $d_1=0.3$ and $d_2=0.01$, and compute $|\Delta \varphi|$, $R$, and $\rho$ to measure PS. The results are presented in Fig. 3.5a, where the blue, orange, and black curves correspond to $|\Delta \varphi|$, $R$, and $\rho$, respectively. For weak levels of $\delta_2$ (e.g. $\delta_2<0.95$ in Fig 3.5a), $\rho$ and $|\Delta \varphi|$ are rapidly decreasing whereas $R$ is rapidly increasing. Then, for strong levels of $\delta_2$ (e.g. $\delta_2>0.95$ in Fig 3.5a), $\rho$ and $|\Delta \varphi|$ are increasing and $R$ is decreasing. That is, PS increasing with the increment of $\delta_2$ for the weak intensity range, reaches an optimal level at an intermediate intensity (e.g. $\delta_2=0.95$ for all measures in Fig. 3.5a), and then begins to decrease as the noise intensity is increased above the optimal intensity value. Such changes in PS can also be observed from the probability density functions of $\Delta \varphi$ shown in Fig. 3.5b, where the dotted blue curve (labelled A), solid blue curve (labelled B), and dashed blue curve (labelled C) correspond to noise intensities $\delta_2=0.05, 0.95, \text{ and } 3$, and points A, B, and C in Fig. 3.5a, respectively. As $\delta_2$ is increased, the peak of the distribution of $\Delta \varphi$ shifts to the right, towards $\Delta\varphi = 0$. For example, when the noise intensity is too weak the peak of the density of $\Delta \varphi$ is slightly less than $\Delta\varphi = 0$, which is marked by curve A in Fig 3.5b. When the noise intensity is increased to the optimal level, the peak of the density of $\Delta \varphi$ moves closer toward $\Delta\varphi = 0$, which is represented by curve B in Fig. 3.5b. And, when the noise intensity is too strong, which is represented by the point and curve C in Fig. 3.5a and 3.5b, respectively, the peak of the density of $\Delta \varphi$ moves farther away from $\Delta\varphi = 0$. Recall, since we are studying the PS of oscillators that are in-phase, values of $\Delta \varphi$ closer to $0$ are indicative of a greater degree PS. Additionally, the distribution of $\Delta \varphi$ is the narrowest for the optimal noise intensity (curve B in Fig. 3.5b) whereas when the intensity of $\delta_2$ is too weak or too strong (i.e. curves A and C) the distribution of $\Delta \varphi$ is considerably wider, which indicates a lesser degree of PS. \par

Thus far, we have considered $\delta_1$ to be fixed at an arbitrary value ($\delta_1 = 0.05$ in Fig. 3.5). To study the effects of $\delta_1$ and $\delta_2$ more systematically, we simulate our model in the parameter space $0.01 \leq \delta_1,\delta_2 \leq 5$ and measure PS using $|\Delta \varphi|$ and $R$, from Eqns. \ref{absolute phase difference} and \ref{R}, respectively. We present the results as heat maps of $|\Delta \varphi|$ and $R$ in Fig. 3.6a and 3.6b, respectively. In Fig. 3.6 warmer colours correspond to larger values and cooler colours correspond to smaller values of $R$ and $|\Delta \varphi|$. In the region of $\delta_1 \in [0.01, 1]$ and $\delta_2 \in [0.1,5]$ we see that PS can be optimized at an optimal (intermediate) level of noise. Indeed, we see that there are multiple routes to PS (i.e. multiple CRs). For example, if we consider a fixed $\delta_1 \in [0.01, 1]$ we observe the characteristic behaviour of CR with the increment of $\delta_2$ in Fig. 3.6a (see the dark blue and bright yellow regions in Fig. 3.6a and 3.6b, respectively). There is an additional region $\delta_1 \in [0.1,5]$ and $\delta_2 \in [0.01,0.05]$ where we see CR. If we consider a fixed $\delta_2 \in [0.01,0.05]$ we see that PS can be maximized at an intermediate value of $\delta_1$. \par
%---------------------------------------------------------------------------------------------------------------
\begin{figure}
\begin{center}
\includegraphics[width=5.8in]{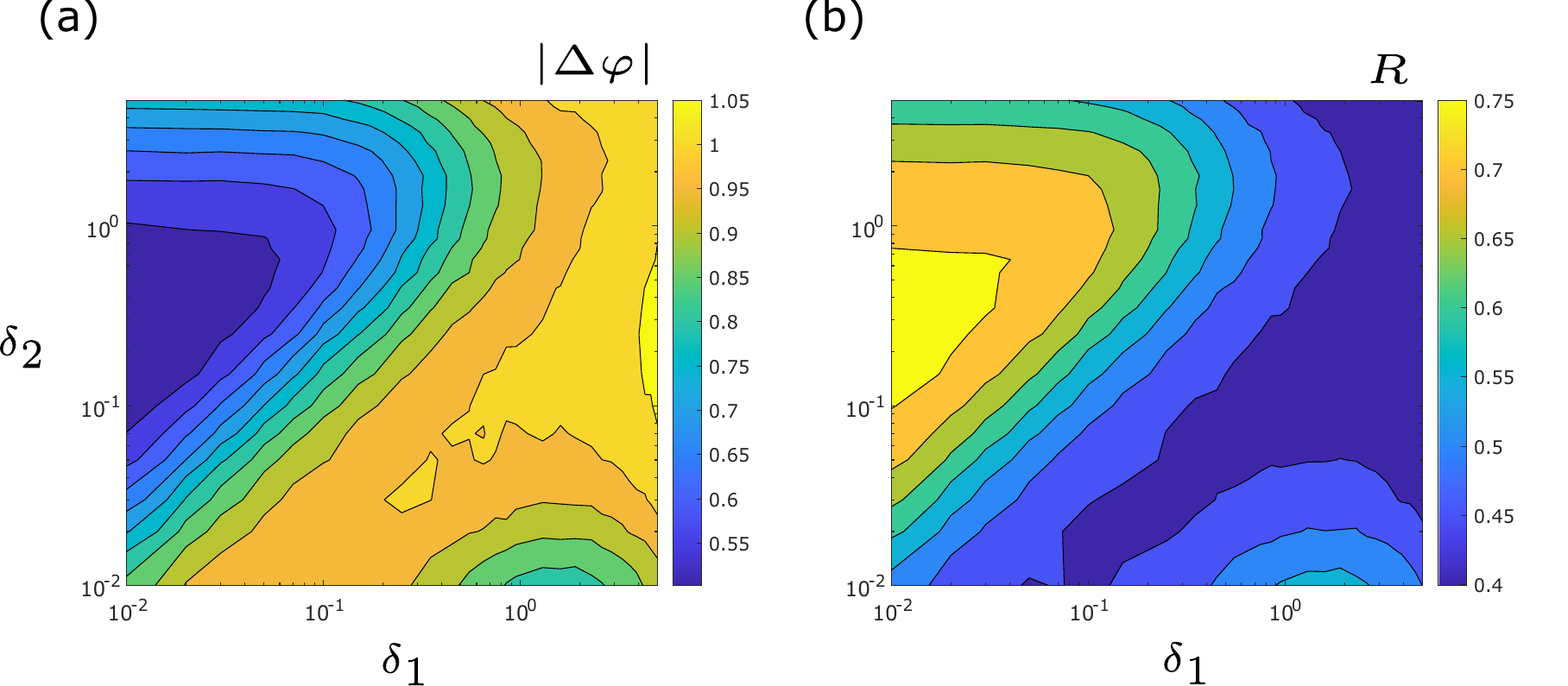}
\caption{Heat maps for (a): absolute CRP difference, $|\Delta \varphi|$, vs. $\delta_1$ vs. $\delta_2$ and (b): mean phase coherence, $R$, vs. $\delta_1$ vs. $\delta_2$. Warmer colours correspond to larger values. Other parameters are: $\alpha=-0.2$; $\gamma=-0.2$; $\omega_0=2$; $\omega_1=0$; $\lambda_0=-0.5$; $d_1=0.3$; and $d_2=0.01$.}
\end{center}
\end{figure}
%----------------------------------------------------------------------------------------------------------------
Recall that the results from Fig. 3.3 in Sec. \ref{noise-induced oscillations} indicated that PS is optimized for asymmetric levels of the intrinsic noise. Considering the ratio $\delta_1/\delta_2$ we see that PS is optimized when $\delta_1/\delta_2 \approx 0.2$ in the region $\delta_1 \in [0.01, 1]$ and $\delta_2\in[0.1,5]$ which agrees with our results in Sec. \ref{noise-induced oscillations}. Moreover, on the line $\delta_1=\delta_2$ (in Fig. 3.6) we see that PS can be increased by increasing $\delta_2$ the region. In other words, if we consider our model with symmetric noise, $\delta_1=\delta_2$, the degree of PS can be increased by changing either $\delta_1$ or $\delta_2$ (i.e. by making the noise asymmetric).

%%%%%%%%%%%%%%%%%%%%%%%%%%%%%%%%%%%%%%%%%%%%%%%%%%%%%%%%%%%%%%%%%%%%%%%%%%%%%%%%%%%%%%%%%%%%%%%%%%%%%%%%%%%%%%%%%%%%%%%%%%%%%%%%%%%%%%%%%%%%%%%%%%%%%%%%%%%%%%%%%%%%%%%%%%%%%%%%%%%%%%%%%%%%%%%%%%%%%%%%%%%%%%%%%%%%%%
\subsection{The effects of $\lambda_0$ on phase synchronization}
\label{sec:lambda}
Recall from Sec. \ref{sec:model} that our model is quiescent and excitable when $\lambda_0<0$. It has been shown that in excitable networks the distance of the control parameter from a critical point (or excitation threshold) has an effect on synchronization (e.g. see \cite{yu2009constructive} and \cite{yu2021noise}). Furthermore, to study the effects of $\lambda_0$ on PS we simulate our model for multiple values of $\lambda_0$ in the excitable regime ($\lambda_0<0$) and study the PS of the system when subject to intrinsic noise. To measure PS we use the absolute CRP difference, $|\Delta \varphi|$, as defined in Eqn. \ref{absolute phase difference}. \par

A series of $|\Delta \varphi|$ curves are shown in Fig. 3.7a, where the solid blue, dashed orange, and dotted black lines correspond to the parameters $\lambda_0=-0.03,-0.5,\text{ and } -1$, respectively. All three $|\Delta \varphi|$ curves show that that PS can be optimized by tuning the intrinsic noise intensity $\delta_2$. At first, the curves in Fig. 3.7a are decreasing, reach a minimum value, and then increase. In other words, all three curves show the characteristic pattern of CR: synchronization is increasing with the increment of $\delta_2$ over a weak intensity range; reaches an optimal point of PS at an intermediate $\delta_2$; and then weakens as $\delta_2$ increases further. Additionally, when $\delta_2$ is in the weak intensity range (e.g. $\delta_2<0.5$ in Fig. 3.7a), $|\Delta \varphi|$ is smaller for values of $\lambda_0$ closer to zero. This tells us that synchronization is enhanced when $\lambda_0$ is closer to the critical point $\lambda_0=0$ (the excitation threshold) over the weak intensity range. The converse is true when the intrinsic noise, $\delta_2$, is strong (e.g. $\delta_2>2$ in Fig. 3.7a), as $\lambda_0$ moves closer to the excitation threshold $\lambda_0=0$, $|\Delta \varphi|$ becomes larger. 

\begin{figure}
\begin{center}
\includegraphics[width=6in]{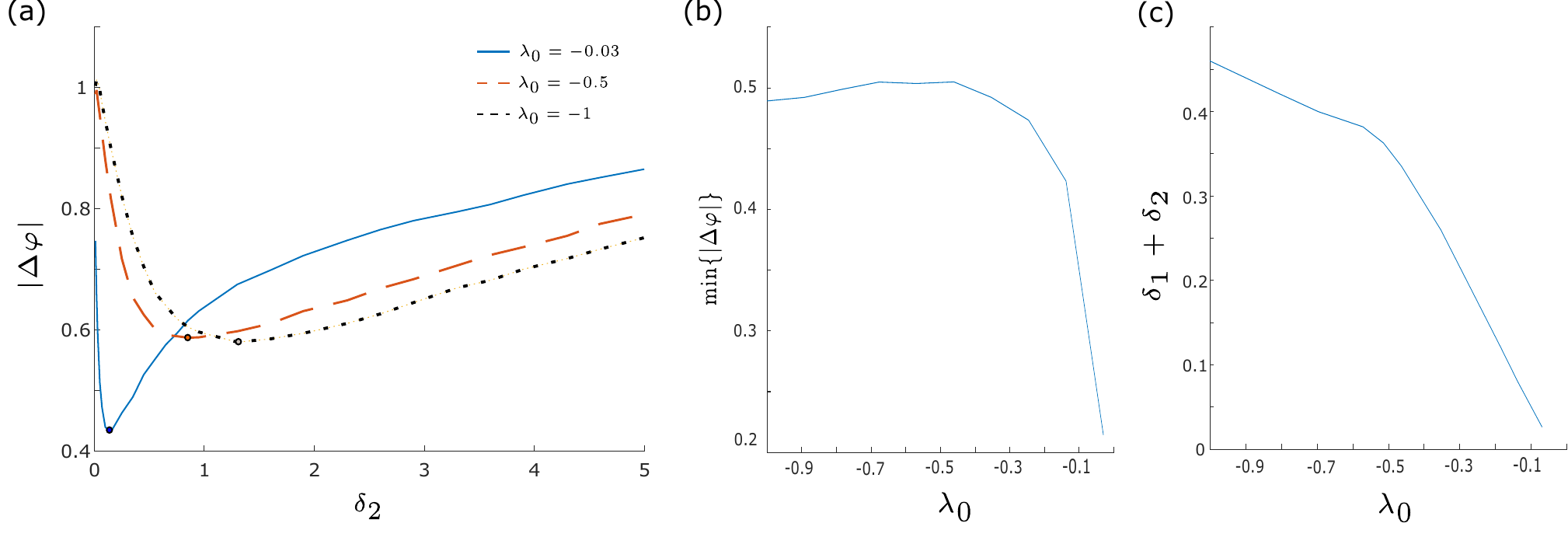}
\caption{(a): Absolute CRP difference, $|\Delta \varphi|$, vs. $\delta_2$ for $\lambda_0=-0.03,-0.5,\text{ and } -1$ with $\delta_1=0.1$. The blue line, orange dashed line, and black dotted line correspond to $\lambda_0=-0.03,-0.5,\text{ and } -1$, respectively. (b): Minimum absolute CRP difference, $\text{min}\{|\Delta \varphi|\}$, vs. $\lambda_0$. (c): Norm-1 of optimal intrinsic noise intensities, $\delta_1+\delta_2$, vs. $\lambda_0$. Note, the Norm-1 is arbitrarily chosen to describe the collective noise strength. Other parameters are: $\alpha=-0.2$; $\gamma=-0.2$; $\omega_0=2$; $\omega_1=0$; $d_1=0.3$; and $d_2=0.01$.}
\end{center}
\end{figure}

From Fig. 3.7a it is evident that the optimal noise intensity and maximum PS are dependent on $\lambda_0$. To better understand this dependence we compute the minimum absolute CRP difference, $\text{min}\{|\Delta \varphi|\}$, and the corresponding optimal intensities $\delta_1$ and $\delta_2$ (note, $\delta_1$ is no longer fixed) for various values of $\lambda_0$ in the excitable regime. The results are presented in Fig. 3.7b and 3.7c which show the $\text{min}\{|\Delta \varphi|\}$ vs. $\lambda_0$ and the norm-1 (or sum) of the corresponding optimal noise intensities, $(\delta_1+\delta_2)$ vs. $\lambda_0$. When $\lambda_0$ is far way from $0$, for example, $-1\leq \lambda_0<-0.5$, the minimum CRP difference and the corresponding optimal intensities are stable. That is, the curves in Fig. 3.7b and 3.7c are relatively flat over this region. This implies that when $\lambda_0$ is sufficiently far from the threshold $\lambda_0=0$, small changes in $\lambda_0$ do not effect the synchronization of our model to a significant degree. When $\lambda_0$ is closer to the excitation threshold, for example, when $-0.5< \lambda_0 \leq0$, in Fig. 3.7b. We observe that as $\lambda_0 \rightarrow 0^-$ the minimum value of $|\Delta \varphi|$ decreases exponentially. This suggests that the PS of our system can be enhanced by shifting $\lambda_0$ closer to the excitation threshold and the optimal synchronization is achieved as $\lambda_0$ is trivially close to $0$ from the left. Lastly, the norm of the optimal $\delta_1$ and $\delta_2$ values decreases as sharply as $\lambda_0 \rightarrow 0^-$. This suggests that as our $\lambda_0$ moves closer to the threshold $\lambda_0=0$ smaller levels of noise are required to produce the most synchronous oscillations.
%-----------------------------------------------------
%%%%%%%%%%%%%%%%%%%%%%%%%%%%%%%%%%%%%%%%%%%%%%%%%%%%%%%%%%%%%%%%%%%%%%%%%%%%%%%%%%%%%%%%%%%%%%%%%%%%%%%%%%%%%%%%%%%%%%%%%%%%%%%%%%%%%%%%%%%%%%%%%%%%%%%%%%%%%%%%%%%%%%%%%%%%%%%%%%%%%%%%%%%%%%%%%%%%%%%%%%%%%%%%%%%%%%
\subsection{The effects of coupling on phase synchronization}
\label{The effects of coupling on synchronization}
%----------------------------------------
%----------------------------------------
\begin{figure}
\begin{center}
\includegraphics[width=6in]{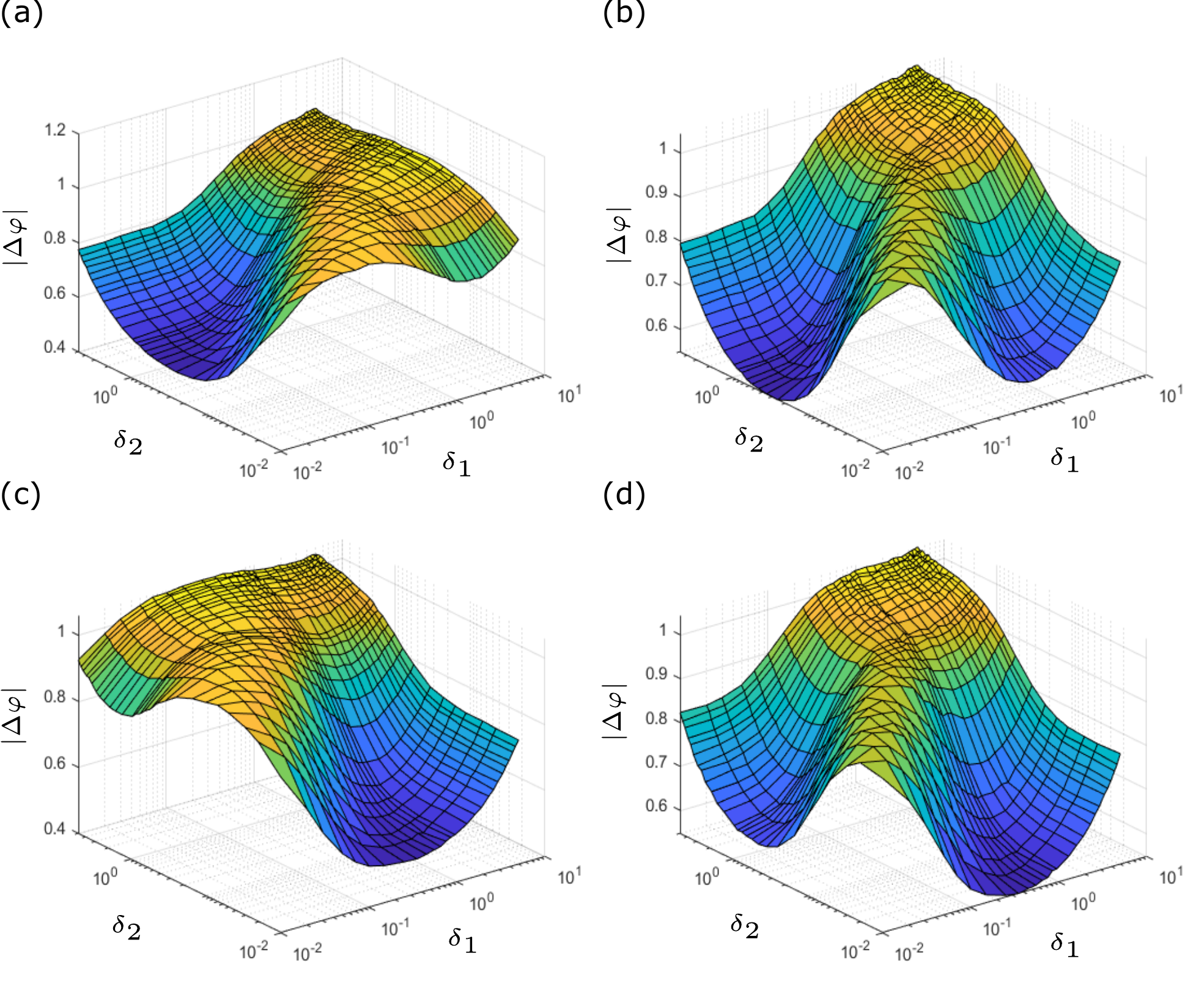}
\caption{Surface plots of the absolute CRP difference, $|\Delta \varphi|$, vs. $\delta_2$ vs. $\delta_1$ for various coupling regimes, $d_1$ and $d_2$. (a): $d_1=0.3$ and $d_2=0.01$. (b): $d_1=0.2$ and $d_2=0.1$. (c): $d_1=0.01$ and $d_2=0.3$. (d): $d_1=0.1$ and $d_2=0.2$. Other parameters are: $\alpha=-0.2$; $\gamma=-0.2$; $\omega_0=2$; $\omega_1=0$; and $\lambda_0=-0.5$.}
\end{center}
\end{figure}
%------------------------------------------
%------------------------------------------
To study the effect(s) of asymmetric coupling on the synchronization of our oscillators we simulate our model with $\lambda_0=-0.5$ under various coupling schemes subject to the intrinsic noise, $\delta_id\eta_i$, $i=1,2$, and measure PS using the absolute CRP difference, $|\Delta \varphi|$, as defined in Eqn. \ref{absolute phase difference}.\par

Fig. 3.8 shows the surface plot of $|\Delta \varphi|$ vs. $\delta_2$ vs. $\delta_1$ for disparate couplings: $d_1>>d_2$ (e.g. $d_1=0.3$ and $d_2=0.01$ in Fig. 3.8a); $d_1>d_2$ (e.g. $d_1=0.2$ and $d_2=0.1$ in Fig. 3.8b); $d_1<<d_2$ (e.g. $d_1=0.01$ and $d_2=0.3$ in Fig. 3.8c); and $d_1<d_2$ (e.g. $d_1=0.1$ and $d_2=0.2$ in Fig. 3.8d). Each panel in Fig. 3.8 displays \textit{two distinct} local minima, which correspond to regions of maximal PS and indicate CR. In these regions, the optimal noise ratios, $\delta_1/\delta_2$, are in the vicinity of $\delta_1/\delta_2 \in [0.025,0.04]$ or $\delta_1/\delta_2 \in [15,35]$. This indicates that PS tends to be maximized when the noise intensities are asymmetric (i.e. $\delta_1/\delta_2 \neq 1$). Indeed, the oscillators exhibit a relativity low degree of PS on the line $\delta_1=\delta_2$ for all panels in Fig. 3.8 (which is in agreement with our results in Sec. \ref{sec:lambda}). Our results further indicate that the ratio $d_1/d_2$ has a significant influence on the ratio $\delta_1/\delta_2$ which optimizes synchronization---the optimal noise ratio. For example, in Fig. 3.8a and 3.8b, when $d_1>d_2$, we see that PS is maximized when $\delta_1<\delta_2$ (or $\delta_1/\delta_2 < 1$) and when $d_1<d_2$ (as in Fig. 3.8c and 3.8d) PS is maximized when $\delta_1>\delta_2$ (or $\delta_1/\delta_2 > 1$) . \par
%------------------------------------------------------------------------------------------------------
\begin{figure}
\begin{center}
\includegraphics[width=6in]{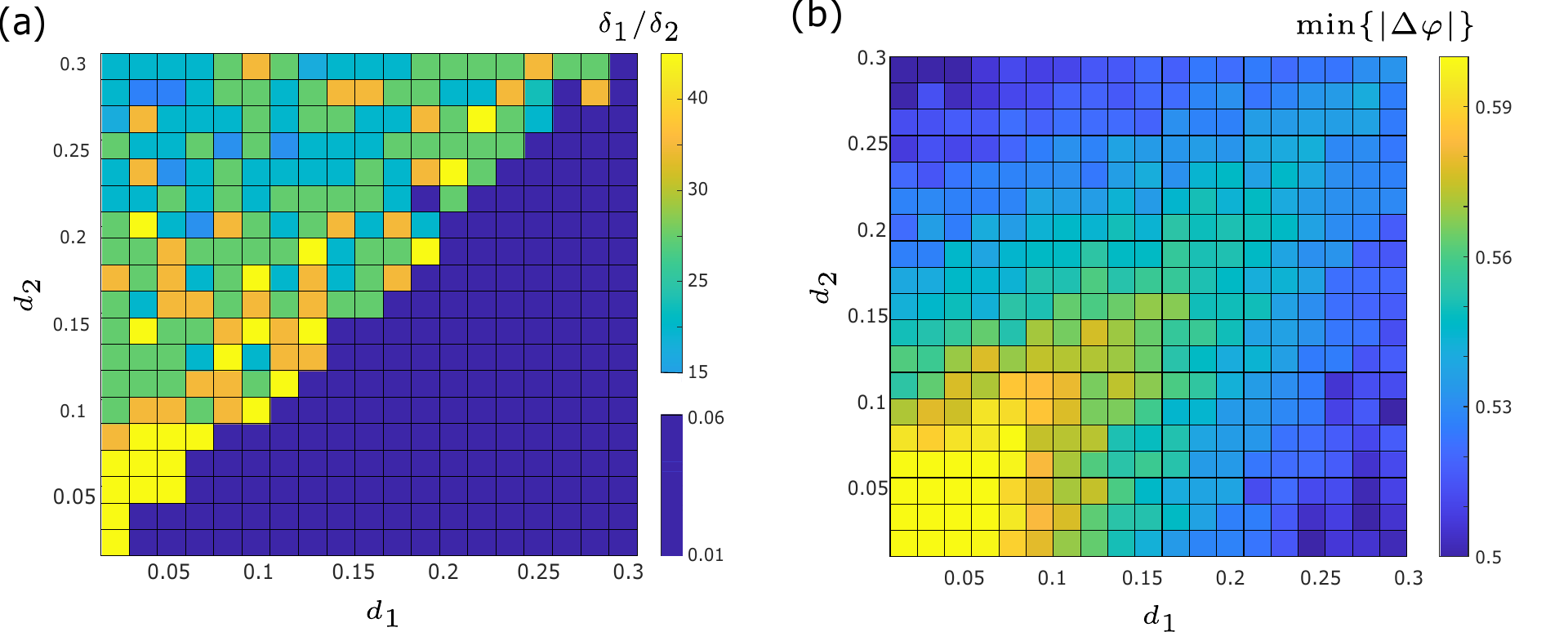}
\caption{Heat maps of  minimum absolute CRP difference $\text{min}\{|\Delta \varphi|\}$ and optimal noise ratio $\delta_1/\delta_2$ vs. coupling strengths, $d_1$ and $d_2$. Panel (a): optimal noise ratio, $\delta_1/\delta_2$, vs. $d_2$ vs. $d_1$ Panel: (b) minimum absolute CRP difference, $\text{min}\{|\Delta \varphi|\}$, vs. $d_2$ vs. $d_1$.  Other parameters are: $\alpha=-0.2$; $\gamma=-0.2$; $\omega_0=2$; $\omega_1=0$; and $\lambda_0=-0.5$.}
\end{center}
\end{figure}
%------------------------------------------------------------------------------------------------------

To explore these results more systematically, we simulate our model for $d_1,d_2 \in [0.01, 0.3]$ and $\delta_1,\delta_2 \in [0.01,0.3]$ and compute minimum absolute CRP difference, $\text{min}\{|\Delta \varphi|\}$ and optimal noise ratio, $\delta_1/\delta_2$, for which the PS is maximized. The results are presented in Fig. 3.9 which displays the heat maps of the optimal noise ratio, $\delta_1/\delta_2$ vs. $d_2$ vs. $d_1$, in panel a, and the heat maps of the minimum absolute CRP difference, $\text{min}\{|\Delta \varphi|\}$ vs. $d_2$ vs. $d_1$, in panel b. We first consider how the choice of coupling affects the optimal noise ratio $\delta_1/\delta_2$. We see from Fig. 3.9a that $\delta_1/\delta_2 \neq 1$ for any considered pair of $d_1$ and $d_2$ . This (again) agrees with our previous results (e.g. Fig. 3.8). However, the implications of Fig. 3.9a are much more robust, as it implies that PS is never maximized when $\delta_1=\delta_2$, independent of the choice of coupling parameters $d_1$ and $d_2$. Furthermore, we see that PS can always be increased by changing either $\delta_1$ or $\delta_2$ when $\delta_1=\delta_2$. Our results demonstrate that there is an interplay between the ratio of coupling strengths and ratio of optimal noise intensities. For example, consider the distinction between the upper and lower triangles of Fig. 3.9a. In the upper triangular region of Fig. 3.9a, where $d_1/d_2<1$ and the optimal noise ratio $\delta_1/\delta_2>1$. Conversely, in the lower triangular region of Fig. 3.9a, where $d_1/d_2>1$ and the optimal noise ratio $\delta_1/\delta_2<1$. Hence, one sees that $d_1/d_2<1 \implies \delta_1/\delta_2>1$ and $d_1/d_2>1 \implies \delta_1/\delta_2<1$. 
Whats-more, our results show that when $d_1/d_2<1$ (or, $d_2>d_1$), the average optimal noise intensities are $\delta_1=0.3462$ and $\delta_2=0.015$ which give an average optimal ratio of $0.3462/0.015 \approx 23.08 $ and when $d_1/d_2>1$ (or, $d_2<d_1$) the average optimal noise intensities are the converse:  $\delta_1=0.015$; and $\delta_2=0.3462$. This gives an average optimal noise ratio of $0.015/0.3465 \approx 0.043$ (note, the latter ratio is the reciprocal the former).\par

Next, we consider how $d_1$ and $d_2$ influence the degree of PS of our model. From Fig. 3.9b we see that when the oscillators symmetrically coupled, PS---as measured by $\text{min}\{|\Delta \varphi|\}$---is positively correlated with the coupling strength $d_1=d_2$. For example, in Fig. 3.9b, when $d_1$ and $d_2$ are symmetric and in the range $[0.01, 0.1]$ (i.e. the yellow region in Fig. 3.9b), $0.56 \leq \text{min}\{|\Delta \varphi|\} \leq 0.65$, whereas, when $d_1$ and $d_2$ are symmetric and in the region $[0.2, 0.3]$ (i.e. top right corner in Fig. 3.9b), $0.5 \leq \text{min}\{|\Delta \varphi|\} \leq 0.53$. Indeed, if we consider the line $d_1=d_2$ on Fig. 3.9b we see that as the coupling strengths increase, the degree of PS increases as well. However, we notice that PS is not maximized near the line $d_1=d_2$. Rather, PS is maximized when the coupling is asymmetric ($d_1 \neq d_2$). In fact, PS is maximized when the ratio $d_1/d_2$ is as large or small as possible. For example, in Fig 3.9b, our model produces the most synchronous oscillations when $d_1/d_2=30 \text{ or } 0.0\overline{3}$. This is in agreement with the results from Fig. 3.8, where we see when $d_1>>d_2$ (Fig. 3.8a) and $d_2>>d_1$ (Fig. 3.8c), the global minimum of $|\Delta \varphi|$ is smaller than when $d_1>d_2$ (Fig. 3.8c) and $d_2>d_1$ (Fig. 3.8d). \par

%%%%%%%%%%%%%%%%%%%%%%%%%%%%%%%%%%%%%%%%%%%%%%%%%%%%%%%%%%%%%%%%%%%%%%%%%%%%%%%%%%%%%%%%%%%%%%%%%%%%%%%%%%%%%%%%%%%%%%%%%%%%%%%%%%%%%%%%%%%%%%%%%%%%%%%%%%%%%%%%%0%%%%%%%%%%%%%%%%%%%%%%%%%%%%%%%%%%%%%%%%%%%%%%%%%%%%%
\section{Discussion}
\label{sec:discussion}
We consider a pair of diffusively coupled $\lambda-\omega$ oscillators with parameters chosen such that the model is in the vicinity of a super-critical HB, quiescent in the absence in noise, and excitable with the addition of an intrinsic noise stimulus. Our results agree with previous studies \cite{ yu2006stochastic,thompson2012stochastic,yu2021noise,yu2008stochastic} which show that noise can play a constructive role in inducing PS in coupled oscillators and that, PS can be optimized by shifting the model closer to the bifurcation point/excitation threshold. The noise-induced PS of excitable systems has been extensively studied in the past (e.g. see \cite{rosenblum2001phase,neiman1999noise,freund2003frequency} and references therein). For simplicity, such studies tend to assume symmetrical interactions between oscillators (i.e. symmetrical coupling and/or stochastic stimulus). However, in biological systems the interactions between the oscillators
are often asymmetric and the assumption of symmetric interactions is strong and likely restrictive (e.g. see \cite{sheeba2009asymmetry,cimponeriu2003inferring}).
\par 
We consider the effect and interplay of asymmetric coupling and asymmetric intrinsic noise on the PS of our model. Our results indicate that PS is maximized when the noise intensities $\delta_1$ and $\delta_2$ are asymmetric (i.e. $\delta_1/\delta_2 \neq 1$). Our results are robust in that we find that the latter is independent of the choice of coupling, $d_1$ and $d_2$. More remarkably, we show that the PS of our model is optimized when the absolute difference between $d_1$ and $d_2$ is as large a possible. We conclude by showing that the asymmetry of coupling and noise are inter-connected such that: $d_1/d_2<1 \implies \delta_1/\delta_2>1$ and $d_1/d_2>1 \implies \delta_1/\delta_2<1$, where $\delta_1/\delta_2$ is the optimal noise ratio and $d_1,d_2>0$.\par

Our study reveals a strong relationship between asymmetric coupling and noise in coupled oscillators and has potential applications in the study of many real-world problems which exhibit asymmetry in interactions between oscillators such as: cardio-respiratory electroencephalogram (EEG) interactions \cite{sheeba2009asymmetry,paluvs2003direction}; optical communication systems and the detection of radar signals in the presence of channel noise \cite{tsang1986bit}; and interactions between ensembles of oscillators in neuronal dynamics \cite{sheeba2009asymmetry,singer1999striving}. Nonetheless, the relationship between the ratios $d_1/d_2$ and $\delta_1/\delta_2$ warrants a deeper investigation. An extension of our current work may be to consider the role of asymmetric noise and coupling in the anti-phase synchronization oscillators of coupled oscillators by considering inhibitory coupling (i.e. $d_i<0$ for an $i \in \{1,2\}$). 

%%%%%%%%%%%%%%%%%%%%%%%%%%%%%%%%%%%%%%%%%%%%%%%%%%%%%%%%%%%%%%%%%%%%%%%%%%%%%%%%%%%%%%%%%%%%%%%%%%%%%%%%%%%%%%%%%%%%%%%%%%%%%%%%%%%%%%%%%%%%%%%%%%%%%%%%%%%%%%%%%%%%%%%%%%%%%%%%%%%%%%%%%%%%%%%%%%%%%%%%%%%%%%%%%%%%%%
%---------------------------------------------------------
%\bibliographystyle{Chicago}

%\bibliography{Bibliography-MM-MC}

%\bibliographystyle{abbrv}
\bibliographystyle{elsarticle-num}
\bibliography{references}

\end{document}